\documentclass[12pt]{article}

\usepackage{color}
\usepackage[table,xcdraw,dvipsnames]{xcolor}
\usepackage{amsmath}
\usepackage{bookmark}
\usepackage{mathrsfs}
\usepackage{amsfonts}
\usepackage{nicefrac}
\usepackage{amssymb}
\usepackage{graphicx,float}
\usepackage{enumitem}
\usepackage{cite}
\usepackage{environ}%

\NewEnviron{neweq}{%
\begingroup
\allowdisplaybreaks
\begin{align}
\begin{split}
    \BODY
    \end{split}
\end{align}
\endgroup}

\NewEnviron{neweq_non}{%
\begingroup
\allowdisplaybreaks
\begin{align*}
    \BODY
\end{align*}
\endgroup}
\usepackage{amsthm}
\newtheorem{theorem}{Theorem}

\def\gavg{\tau_{\rm m}}

\def\e{{\rm e}}
\def\re{{\rm Re}}
\def\im{{\rm Im}}
\def\unig{{g_{\rm u}}}
\def\gammag{{g^{\rm \gamma}_{p}}}

\def\ds{{\rm d}s}

\newcommand\ddc[2]{\cfrac{{\rm d} #1}{{\rm d} #2}}

\definecolor{Nblue}{HTML}{0080ff}
\definecolor{NewGreen}{rgb}{0, 0.501, 0}
\definecolor{Red1}{rgb}{0.858, 0.188, 0.478}
\usepackage{hyperref}
\hypersetup{
    colorlinks=true,
    citecolor=Nblue,
    linkcolor=NewGreen,
    filecolor=magenta,      
    urlcolor=Red1
    }

\usepackage
[left=2cm,right=2cm,top=2cm,bottom=3cm]
{geometry}

\makeatletter
\def\myrulefill{\leavevmode\leaders\hrule height .7ex width 1ex depth -0.6ex\hfill\kern\z@}
\makeatother

\usepackage{silence}
\date{}
\title{Stability Analysis  of A Single-Species Model with Distributed Delay}
\author{Isam Al-Darabsah
\footnote{Department of Mathematics and Statistics, Faculty of Science and Arts, Jordan University of Science and Technology, P.O. Box 3030, Irbid 22110, Jordan}
	\textsuperscript{,}\footnote{Corresponding author. Email: imaldarabsah@just.edu.jo, ialdarabsah@gmail.com}}
\begin{document}

\maketitle

\begin{abstract}
The logistic equation  has many applications and is used frequently in different fields, such as biology, medicine, and economics.
In this paper, we study the stability of a single-species logistic model with a general distribution delay kernel and an inflow of nutritional resources at a constant rate.
In particular, we provide precise conditions for the linear stability of the positive equilibrium and the occurrence of Hopf bifurcation. 
We apply the results to three delay distribution kernels: 
Uniform, Dirac-delta, and gamma distributions. 
Without an inflow, we show that the positive equilibrium is stable for a relatively small delay and then loses its stability through the Hopf bifurcation when the mean delay $\gavg$ increases with the three distributions.
In the presence of an inflow, the model dynamics depend on the delay distribution kernel.
In the uniform and Dirac-delta distributions cases, we find that the dynamics are similar to the absence of a nutrient influx.
In contrast, the dynamics depend on the delay order $p$ when considering the gamma distribution.
For $p=1$, the positive equilibrium is always stable. 
While for $p=2$ and $p=3$, we find stability switching of the positive equilibrium resulting from the increase of the value of $\gavg$, where the positive equilibrium is stable for a relatively short period; then, it loses stability via Hopf bifurcation as $\gavg$ increases; after then, it stabilizes again with an increase in $\gavg$. 
The main difference between the delay orders $p=2$ and $p=3$ is that for relatively large $\gavg$ and intrinsic growth rate, the positive equilibrium  can be stable   when $p=2$, but it will be unstable when $p=3$.

\bigskip

\noindent \textbf{keywords}: { Single-species model $\cdot$  Stability   $\cdot$  Distributed delay}

\end{abstract}

	\section{Introduction}
	\label{sec_introduction}

Time delay can be incorporated into  ecological models to represent a time lag in different biological processes.
For instance, 
time lag due to maturity period \cite{al2019dynamic}, incubation time \cite{liz2014delayed},
reproductive process time \cite{wangersky1956time}, 
and  reaction time of predation \cite{dubey2019global}.
Time delays  can take different forms:  
 Fixed time delay is used when the time lag is the same for all population members as in maturation time \cite{al2019dynamic}.
To consider the variation among the population members, the distributed delay allows for a more appropriate description of the time lag than a fixed delay to reflect that the time lag is not the same for everyone in the population, but  
vary according to a distribution \cite{campbell2009approximating,lin2018alternative}.
The distributed delay kernel can take different forms, such as  the Dirac-Delta function, uniform distribution, and  gamma distribution. 
Another form is a time-dependent delay, which is used when the time lag is influenced by certain factors, such as temperature 
 \cite{bartuccelli1997population,al2018periodic}.
Other researchers used state-dependent delay to study the maturation of a stem cell population \cite{getto2016differential} and stochastic delay to analyze gene regulatory networks \cite{gomez2016stability}.
In general, introducing time delay in ecological models exhibits complex dynamics compared to ordinary differential equations, as it can destabilize equilibrium points and give rise to the appearance of limit cycles \cite{balachandran2009delay,zhao2003dynamical}.


In $1838$, Verhulst  introduced the classical logistic ordinary differential equation (ODE)  model to describe population growth in a limited environment.  
After Pearl and Reed rediscovered it in the $1920$s \cite{pearl1920rate}, it became a valuable tool in mathematical ecology, where it was applied to model population dynamics such as bacteria, cells, and human or animal populations with limited nutrients. 
In $1948$, Hutchinson \cite{hutchinson1948circular} modified Verhulst's  model to a logistic delay differential equation (DDE)  and incorporated a fixed delay in the density-dependent feedback on population dynamics.
The incorporated delay represented the time lag between the instant when the population reaches a certain level and the moment when the effective reproductive rate is updated. 
It has been shown that when the time delay increased for large values, an oscillation arose in population density via Hopf bifurcation \cite{beretta1987global,ruan2006delay}.


Different single-species logistic models have been built by incorporating additional biological processes into the classical Verhulst model \cite{ruan2006delay,song2006stability,li2002periodic,li2019stability,zhang2013single,liu2013note,liu2016analysis,sawada2022stability,tarasov2019logistic}. For example, 
in \cite{ruan2006delay},   the dynamics of logistic equation with fixed time delay was studied. 
The author found that a large time delay can cause the positive equilibrium to become unstable and lead to the formation of a stable limit cycle.
In \cite{li2002periodic}, the authors investigated how the seasonality of the changing environment can impact population growth by considering a state-dependent delay logistic equation. 
In \cite{li2019stability}, a single-species model with a constant  harvesting rate and weak delay kernel was considered. 
The authors found that 
In \cite{liu2013note}, the authors discussed the stability of 
stochastic logistic equation. They demonstrated that the stability of the positive equilibrium is negatively impacted by noise.
In \cite{tarasov2019logistic}, the authors applied a logistic equation with distributed lag on  economics. 
Recently, in \cite{sawada2022stability},  the authors studied the stability of the positive equilibrium of the logistic model with a gamma distribution kernel.   
When the delay order in the gamma distribution is two, they showed that the positive equilibrium changes to be unstable from being stable first and returning to being stable again  through Hopf bifurcation  by increasing the mean time delay.   
However, when the delay order is three, the positive equilibrium losses its stability and becomes unstable for large mean time delay.
In this work, we consider a general form of delay distribution and study the stability and the occurrence of Hopf bifurcation. 
We also apply the results to three delay distribution kernels: Uniform, Dirac-delta, and gamma distributions.
Then, we compare The results  with the parts in the literature.


The paper is organized as follows. In Section \ref{sec_model}, we provide the mathematical model and the stability analysis of the positive equilibrium point. 
In Section \ref{sec_applications}, we provide three distributions, uniform, Dirac-delta, and gamma, and discuss their biological meaning. 
Then, we apply stability results to all three distributions and compare results to what existed in the literature.
We discuss our results in Section \ref{sec_conclusions}.

\section{Mathematical model and stability analysis}
\label{sec_model}

Let $n(t)$ be the density of an organism s population at time $t$. Assume that the organism grows in an environment with constant nutritional resources.
Then, the dynamics of $n(t)$ can be represented by 
a single-species model with logistic growth  \cite{ruan2006delay} 
\begin{equation} \label{model_v0}
    \ddc{n(t)}{t}=r\, n(t)\, \left[   
1-\cfrac{1}{K}\int_{-\infty}^{t} n(s)\, g(t-s)\,\ds
    \right].
\end{equation}
where $r$ is the intrinsic growth rate, and $K$ is the carrying  capacity of the resources.
The function   $g(\cdot)$ is the  kernel of the delay distribution with compact support, that is,  
\[g(s)\ge 0\qquad\text{and}\qquad \int _{0}^{\infty }g(s)\ds=1.\]
We calculate the \textit{mean delay} as
\begin{align*}
     \gavg=\int_0^{\infty} s g(s) ds.
 \end{align*}
The function $g$ states that the population growth will be proportionate to the size of the population in the past and will solely depend on individuals who can survive the delay. 

\noindent Assume there is an inflow of more nutritional resources at a constant rate of $D$. Then, model \eqref{model_v0} has the form \cite{sawada2022stability}
\begin{equation} \label{model_v1}
    \ddc{n(t)}{t}=r\, n(t)\, \left[   
1-\cfrac{1}{K}\int_{-\infty}^{t} n(s)\, g(t-s)\,\ds
    \right]+D.
\end{equation}

\noindent  By setting $\bar{s}= t-s$, model \eqref{model_v1} can be written, after dropping the bars,  as 
\begin{equation} \label{model}
    \ddc{n(t)}{t}=r\, n(t)\, \left[   
1-\cfrac{1}{K}\int_{0}^{\infty} n(t-s)\, g(s)\,\ds
    \right]+D.
\end{equation}

In the rest of the manuscripts, we study the stability and the existence of Hopf bifurcation of the model   \eqref{model}  with a general distribution kernel.

When $D\ge0$, model \eqref{model} has  only one positive equilibrium 
\[n^*=\cfrac{\left(1+\sqrt{1+\cfrac{4D}{rK}}\right)K}{2}.\]
Notice that when  $D=0$, then  $n^*=K$. Moreover, in this case, the  trivial equilibrium $n=0$ exists, and it is unstable due to the positive eigenvalue $\lambda=r$.

Define $\tilde{n}=n-n^*$. Then the linearization of \eqref{model} at $n^*$ is
\begin{equation}\label{linear_eq}
 \ddc{\tilde{n}(t)}{t}= r\left( \cfrac{K-n^*}{K}\right)\tilde{n}(t)
-\cfrac{rn^*}{K} \int_0^{\infty}\tilde{n}(t-s)\,g(s)\,\ds. 
\end{equation}

For the zero delay case $\gavg=0$, i.e., ${g}(s)=\delta_{0}(s)$,  system \eqref{linear_eq} becomes
     \[\ddc{\tilde{n}(t)}{t}=-r\left( \cfrac{2n^*-K}{K}\right)\tilde{n}(t)\]
and the characteristic equation  is $\lambda+ r(2n^*/K-1)=0$. Hence, 
 the  equilibrium $n^*$ is locally asymptotically stable  
due to $0<K\le n^*$ and  $D\ge0$.

Let $\gavg>0$. 
To study the dependence of the linear stability of $n^*$ on the mean delay $\gavg$, we rescale  
 the dimensional variables  as  $\bar{t}=t/\gavg$ and $\bar{s}=s/\gavg$. 
After dropping the bars, model \eqref{linear_eq} becomes
\begin{equation} \label{linear_eq_non_bim}
  \ddc{\tilde{n}(t)}{t}= \gavg r\left( \cfrac{K-n^*}{K}\right)\tilde{n}(t)
- \gavg \cfrac{rn^*}{K} \int_0^{\infty}\tilde{n}(t-s)\,\widehat{g}(s)\,\ds.
\end{equation}
where $\widehat{g}(s)=\gavg\, g(\gavg s)$.
By applying the Laplace transform to \eqref{linear_eq_non_bim} (with zero initial condition),  the characteristic equation can be written as
\begin{equation}\label{characteristic_eq}
    \Delta(\lambda)=\lambda+\gavg r\left( \cfrac{n^*-K}{K}\right)+\gavg \cfrac{rn^*}{K}\, \widehat{G}(\lambda)=0
\end{equation}
where
\[\widehat{G}(\lambda)=\int_0^{\infty}\e^{-\lambda s}\,\widehat{g}(s)\,\ds\] 
is the Laplace transform of $\widehat{g}$.

Recall that $n^*$ is locally asymptotically stable when $\gavg= 0$. 
We seek conditions on $\gavg$ such that $\re(\lambda)$ changes its sign as $\gavg$ increases. 
In other words, the characteristic equation
\eqref{characteristic_eq} must have a pair of pure imaginary eigenvalues.

It is clear that $\lambda=0$ is not an eigenvalue value because  $\Delta(0)=r(2n^*/K-1)>0$. 
To determine the existence of a pair of pure imaginary  eigenvalues, we substitute $\lambda= i\omega$  ($\omega>0$ and $i=\sqrt{-1}$)   in $\Delta(\lambda)$ defined in \eqref{characteristic_eq}. Consequently, we have
\[\Delta(i\omega)=i\omega+\gavg r\left( \cfrac{n^*-K}{K}\right)+\gavg \cfrac{rn^*}{K}\, \int_0^{\infty}\e^{-i\omega s}\,\,\widehat{g}(s)\,\ds=0\]
Separating the real and imaginary parts results in 
\begingroup
\allowdisplaybreaks
\begin{equation}\label{Re_Im}
  \cfrac{K-n^*}{K}  =   \cfrac{ n^*}{K}\, C(\omega) 
 \quad \text{and}\quad
 \omega=  \gavg \cfrac{rn^*}{K}\, S(\omega),
\end{equation}
\endgroup
where $\widehat{G}(i\omega)=C(i\omega)-i\,S(i\omega)$ with
\begin{equation}\label{CSSS}
C(\omega):=\int_{0}^{\infty}\cos \! \left(\omega  s \right) \widehat{g} \! \left(s \right)\ds\quad \text{and}\quad
S(\omega):= \int_{0}^{\infty}\sin \! \left(\omega  s \right) \widehat{g} \! \left(s \right)\ds.
\end{equation}
Now we study the  influence of varying the mean delay $\gavg$ on the stability of the positive equilibrium $n^*$.

 \subsection[The case of D=0]{The case of $D=0$.}
Recall that when $D=0$, the positive equilibrium $n^*=K$.  Hence, the characteristic equation \eqref{characteristic_eq} becomes 
\begin{equation}\label{characteristic_eq_D0}
 \Delta(\lambda)=\lambda+\gavg r\, \widehat{G}(\lambda)=0   
\end{equation}
and 
 equation \eqref{Re_Im} implies 
 \begin{equation}\label{eq_CS_1}
  C(\omega)=0 \quad \text{and}\quad S(\omega)= \cfrac{\omega}{\gavg r}.
 \end{equation}
To determine how the sign of eigenvalues changes as $\gavg$ increases, 
we calculate the rate of change of the real part of $\lambda$ with respect to $\gavg$. To this end, firstly, notice that 
\[\widehat{G}'(\lambda)=-\int_0^{\infty}s\e^{-\lambda s}\,\widehat{g}(s)\,\ds\] 
Hence,
\[\widehat{G}'(i\,\omega)=-S'(\omega)-i\, C'(\omega).\] 
From  \eqref{characteristic_eq_D0}, we obtain that 
\[1+r \left(\gavg \widehat{G}'(\lambda) + \widehat{G}(\lambda) \ddc{\gavg}{\lambda}\right)\, =0.\]  
Now we have 
\begin{neweq_non}
 \left.\left(\ddc{\lambda}{\gavg} \right)^{-1}  \right|_{\lambda=i\omega}&=
-\left. \left( \cfrac{1+r\, \gavg\,  \widehat{G}'(\lambda)}{r\,  \widehat{G}(\lambda)} \right) \right|_{\lambda=i\omega}\\
&=-  \cfrac{1-r\, \gavg\, (S'(\omega)+i\, C'(\omega)) }{r\, (C(\omega)-i\, S(\omega))}     
\end{neweq_non}
Using \eqref{eq_CS_1}, we get  
\begin{neweq_non}
\left.\re\left(\ddc{\lambda}{\gavg} \right)^{-1}  \right|_{\lambda=i\omega}
&=-\cfrac{  C'(\omega)}{r\, S(\omega)}
=- \cfrac{\gavg}{\omega}C'(\omega).
\end{neweq_non}
Thus crossing the imaginary axis through the solution of \eqref{eq_CS_1} depends on the sign of $C'(\omega)$.   
Therefore,  when crossing the imaginary  axis, $\re(\lambda)$ changes from negative to positive (resp. positive to negative) when $C'(\omega)<0$ (resp. $C'(\omega)>0$ ). Hence, there exists $\gavg^*>0$ such that  a Hopf bifurcation occurs at $\gavg^*$.

Recall that  $n^*=K$ is   locally asymptotically stable when $\gavg=0$. Consequently, we have the following result.
\begin{theorem}\label{thm_d_0}
    Let $\omega_0$ be the smallest  solution of \eqref{eq_CS_1} such that   $C'(\omega_0)<0$. Then,  a Hopf bifurcation occurs at $ \gavg^*=\omega_0/ (r\,S(\omega_0))$. 
    Consequently, $n^*=K$ is locally asymptotically stable for    
\begin{equation}\label{stability_cond_1}
    0< \gavg <\cfrac{\omega_0}{r\, S(\omega_0)}  
\end{equation}
 and it is unstable when $\gavg >\frac{\omega_0}{r\, S(\omega_0)}$. 
\end{theorem}

 \subsection[The case of D>0]{The case of $D>0$.} 
Since $K<n^*$, equation \eqref{Re_Im} gives  
 \begin{equation}\label{eq_CS_2}
-1<C(\omega)= \cfrac{K}{n^*}-1 <0\quad \text{and}\quad  S(\omega)=\cfrac{\omega K}{\gavg r\, n^*}.
 \end{equation}

Now we study the  rate of change of $\re(\lambda)$ with respect to $\gavg$.
Firstly, from \eqref{eq_CS_2} we obtain
\begin{neweq_non}
\gavg=-\cfrac{K\, \omega}{r(n^*-K)}\, \, \cfrac{C(\omega)}{S(\omega)}.
\end{neweq_non}
Hence,
\begin{neweq_non}
\ddc{\gavg}{\omega}=-\cfrac{K}{r(n^*-K)}\, \, \cfrac{1}{S(\omega)} 
\left( C(\omega)+\omega\cfrac{C'(\omega)\, S(\omega)-C(\omega)S'(\omega)}{S(\omega)}    \right).
\end{neweq_non}
Using $rS(\omega)/K= \omega/\gavg\, n^*$ in \eqref{eq_CS_2}, we have 
\begin{neweq_non}
\ddc{\gavg}{\omega}=-\cfrac{\gavg\,n^*}{(n^*-K)\omega}\,  
\left( C(\omega)+\omega\cfrac{C'(\omega)\, S(\omega)-C(\omega)S'(\omega)}{S(\omega)}    \right).
\end{neweq_non}
Form the characteristic equation \eqref{characteristic_eq}, we have 
\[1+\gavg r\left( \cfrac{n^*-K}{K}\right)+ \cfrac{n^*}{K} \left(\gavg \widehat{G}'(\lambda) + \widehat{G}(\lambda) \ddc{\gavg}{\lambda}\right)\, =0.\]  
Consequently, we obtain 
\begin{neweq_non}
\left.\re\left(\ddc{\lambda}{\gavg} \right)^{-1}  \right|_{\lambda=i\omega}
&=-  \,  
\re\left(\cfrac{K+rn^*\gavg\, \widehat{G}'(\lambda)}{r(n^*-K)+r\,n^*\,\widehat{G}(\lambda)}\right)\\
&=-  \,  
\re\left(\cfrac{K-rn^*\gavg\,  (S'(\lambda)+i\,C'(\lambda))}{r(n^*-K)+r\,n^*\,\,  (C(\lambda)-i\,S(\lambda))}\right)\\
&=-\cfrac{(n^*-K)(K-r\gavg n^*S'(\omega))}{q^2(\omega)}\\
&\qquad-\cfrac{Kn^*(C(\omega)+r\gavg n^*/K(C'(\omega)S(\omega)-C(\omega)S'(\omega)))}{q^2(\omega)}
\end{neweq_non}
where $q^2(\omega)= (r(n^*-K)+rn^*C(\omega))^2+r^2{n^*}^2S^2(\omega) $. 
Using ${\gavg r\, n^*}/{K}=\omega/S(\omega)$ in \eqref{eq_CS_2}, we obtain
\begin{neweq}\label{re_lambda}
\left.\re\left(\ddc{\lambda}{\gavg} \right)^{-1}  \right|_{\lambda=i\omega}
&=\cfrac{r(n^*-K)}{q^2(\omega)}\left( r\gavg n^*S'(\omega)-K+\cfrac{K\omega}{\gavg}\,\ddc{\gavg}{\omega} \right).
\end{neweq}
Thus crossing the imaginary axis through the solution of \eqref{eq_CS_1} depends on the sign of $S'(\omega)$ and $\ddc{\gavg}{\omega}$.   
Therefore,  if
\begin{equation}\label{condition_d_ge_1}
   r\gavg n^*S'(\omega)+\cfrac{K\omega}{\gavg}\,\ddc{\gavg}{\omega}
 \begin{array}{cc} > \\ < \\  \end{array}
K, 
\end{equation}
then $\re(\lambda)$ crosses the imaginary  axis from left to right (resp. right to left).  
Notice that we require knowledge of the distribution $\widehat{g}(s)$ to obtain an explicit condition for how the
  eigenvalues change when \eqref{eq_CS_2} holds.

\begin{figure*}[hbt!]
     \centering
         \includegraphics[width=1\textwidth]{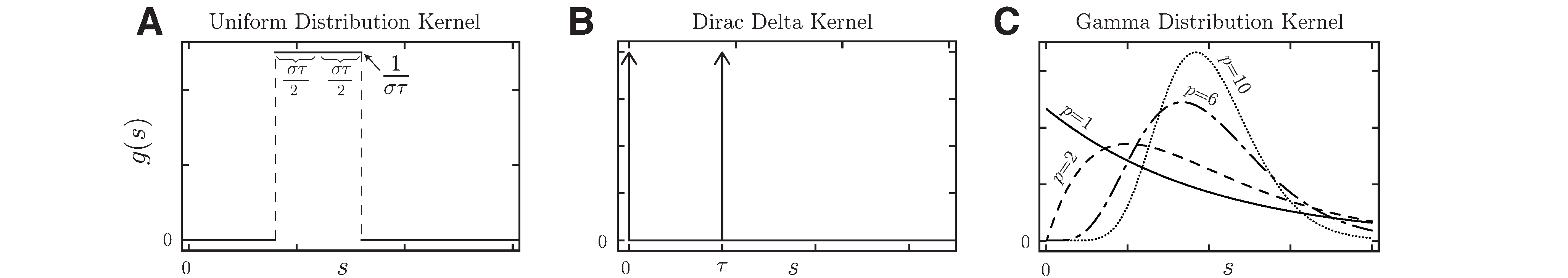}
          \caption{\textbf{Delay distribution kernel.}  
        \textbf{(A)} Uniform distribution kernel $\unig(s)$ in \eqref{kernel_uni}.
        \textbf{(B)} Dirac-Delta kernel $\delta_{0}(s)$ and $\delta_{\tau}(s)$ in \eqref{kernel_delta}.
        \textbf{(C)} Gamma distribution kernel $\gammag(s)$ when $p=1,2,6$, and $10$  in \eqref{kernel_gamma}.
}
  \label{Fig:kernel_forms}
\end{figure*}

\section{Applications}
\label{sec_applications}

 The delay distribution kernel $g(s)$  can take different forms. In this section, we apply the results in Section \ref{sec_model} to different distribution kernels.

 \subsection{Application 1: Uniform distribution kernel}	
\label{sec_uniform}

The uniform distribution kernel (Fig. \ref{Fig:kernel_forms}A) can be written as: 
   \begin{equation}\label{kernel_uni}
         g(s)=\unig(s):=\left\{\begin{array}{lll}\cfrac{1}{\sigma\tau}& \text{if}& \tau (1 -\frac{\sigma}{2}) \le s\le \tau (1 + \frac{\sigma}{2}), \\ 
         &&\\
         0& \text{if} & \text{otherwise}. \end{array}\right.
     \end{equation}
The parameter $\sigma\in(0,2)$ controls the width and height of the distribution with 
the mean time delay  $\gavg=\tau$.
In this case, model \eqref{model} reduces to an integro-differential equation (IDE) of the form
\begin{equation} \label{model_uni}
    \ddc{n(t)}{t}=r\, n(t)\, \left[   
1-\cfrac{1}{\sigma\tau K}\int_{\tau(1-\sigma/2)}^{\tau(1+\sigma/2)} n(t-s)\,\ds
    \right]+D.
\end{equation}
From a biological point of view, the distribution $\unig(s)$  means that the maximum influence on the population density  at the present time $t$  depends equally likely on the population density at any previous time $t-s$.

The normalized uniform distribution has the form 
 \begin{equation}\label{kernel_uni_nor}
         \widehat{g}(s)=\widehat{g}_{\rm u}(s):=\left\{\begin{array}{lll}\cfrac{1}{\sigma}& \text{if}&  (1 -\frac{\sigma}{2}) \le s\le  (1 + \frac{\sigma}{2}), \\ 
         &&\\
         0& \text{if} & \text{otherwise}. \end{array}\right.
     \end{equation}
where $0<\sigma< 2$. 
Then, the linearization equation around $n^*$ can be written as 
\begin{equation} \label{linear_eq_non_bim_uni}
  \ddc{\tilde{n}(t)}{t}= \gavg r \left( \cfrac{K-n^*}{K}\right)\tilde{n}(t)
-  \cfrac{rn^*\gavg}{K\sigma} \int_{1-\sigma/2}^{1+\sigma/2}\tilde{n}(t-s)\,\ds.
\end{equation}
Recall that $\gavg=\tau$ in the uniform distribution kernel. Consequently, the characteristic equation is
 \begin{equation}\label{characteristic_eq_uni}
    \Delta(\lambda)=\lambda+\tau r\left( \cfrac{n^*-K}{K}\right)+\cfrac{2\, r\, n^*\, \tau}{K\, \sigma} \cfrac{\sinh(\sigma\lambda/2)\, \e^{-\lambda}}{ \lambda}=0
\end{equation}
with 
\[
C(\omega)=\cfrac{1}{\sigma}\int_{1-\sigma/2}^{1+\sigma/2}\cos(\omega s)\ds=\cfrac{2\cos(\omega)\sin(\sigma\omega/2)}{\sigma\omega}\]
and
\[S(\omega)=\cfrac{1}{\sigma}\int_{1-\sigma/2}^{1+\sigma/2}\sin(\omega s)\ds=\cfrac{2\sin(\omega)\sin(\sigma\omega/2)}{\sigma\omega}.\]

\subsubsection[The case of D=0]{The case of $D=0$.} In this case, the curves of pure imaginary eigenvalues are 
\begingroup
\allowdisplaybreaks
\begin{equation}\label{Re_Im_uni_1}
\cfrac{2\cos(\omega)\sin(\sigma\omega/2)}{\sigma\omega}
=0 
 \quad \text{and}\quad
  \cfrac{2\sin(\omega)\sin(\sigma\omega/2)}{\sigma\omega}= \cfrac{\omega }{\tau r}.
\end{equation}
\endgroup
From the first equation of \eqref{Re_Im_uni_1}, we have
\[\cos(\omega)=0 \,\,\, \text{or}\,\,\, \sin(\sigma\omega/2)=0.\]
Thus
\[\omega=\cfrac{\pi}{2}+k_1\pi
\,\,\, \text{or}\,\,\,
 \omega=\cfrac{2k_2\pi}{\sigma},\qquad k_1,k_2=0,1,2,\ldots. \]
Hence, the smallest positive root is $\omega_0=\pi/2$ due to   $\sigma\in(0,2)$.
Notice that
\[C'(\pi/2)=-\cfrac{4\sin(\pi\sigma/4)}{\pi\sigma}<0.\]
Thus, it follows by  Theorem \ref{thm_d_0} that $\re(\lambda)$ changes from negative to positive when crossing the imaginary axis, and a Hopf bifurcation occurs at $\tau={\pi^2\sigma }/{(8\,r\,   \sin(\sigma\pi/4))}$. 
Consequently, from \eqref{stability_cond_1}, we know that  the equilibrium $n^*=K$ is   locally asymptotically stable when 
\[0<\tau  < \cfrac{\pi^2\sigma }{8\, r\, \sin(\sigma\pi/4)}\]
and unstable when 
\[ \tau > \cfrac{\pi^2\sigma }{8\, r\, \sin(\sigma\pi/4)}.\]


 \subsubsection[The case of D>0]{The case of $D>0$.} 
 The curves of pure imaginary eigenvalues are 
 \begingroup
\allowdisplaybreaks
\begin{equation*} 
-1<\cfrac{2\cos(\omega)\sin(\sigma\omega/2)}{\sigma\omega}
=\cfrac{K}{n^*}-1 <0 
 \qquad\text{and}\qquad 
  \cfrac{2\sin(\omega)\sin(\sigma\omega/2)}{\sigma\omega}= \cfrac{\omega K}{\tau r\, n^*}.
\end{equation*}
\endgroup
Dividing the two equations gives the equation 
\begin{equation}\label{tan_w}
    \tan(\omega)= -\cfrac{K(n^*-K)}{\tau r}\, \omega. 
\end{equation}
Consequently, considering $\omega$ as a parameter in $(\frac{\pi}{2},\pi)$ or $(\frac{3\pi}{2},2\pi)$. 
Since $\tau$ is positive and $K/n^*-1$ is negative, the only part of interest is the one lying in the first quadrant, that is, $\omega\in(\frac{\pi}{2},\pi)$.

\begin{figure*}[hbt!]
     \centering
         \includegraphics[width=1\textwidth]{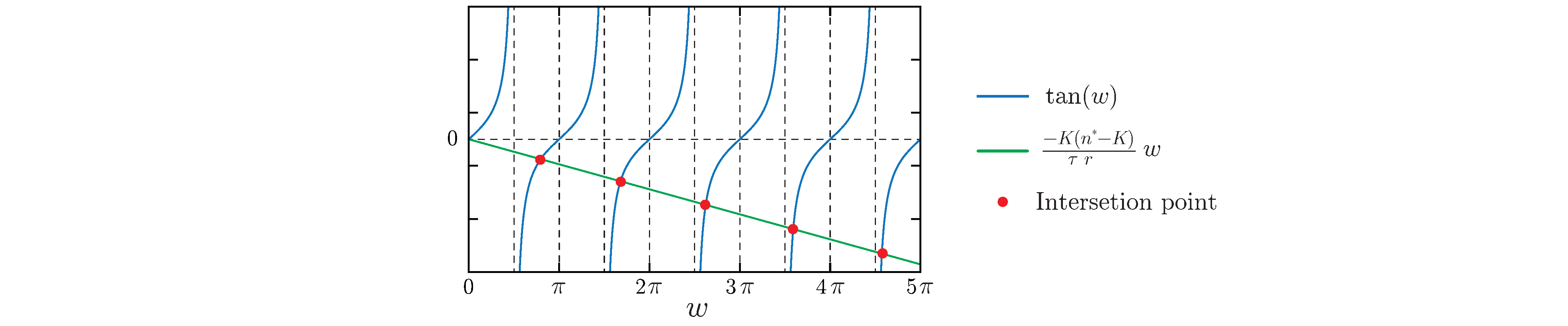}
         \caption{\textbf{Illustration of the existence of pure imaginary eigenvalues when $D>0$.}  
          Roots in equation \eqref{tan_w}.}
  \label{Fig4_tan_w}
\end{figure*}

\begin{figure*}[hbt!] 
     \centering
         \includegraphics[width=1\textwidth]{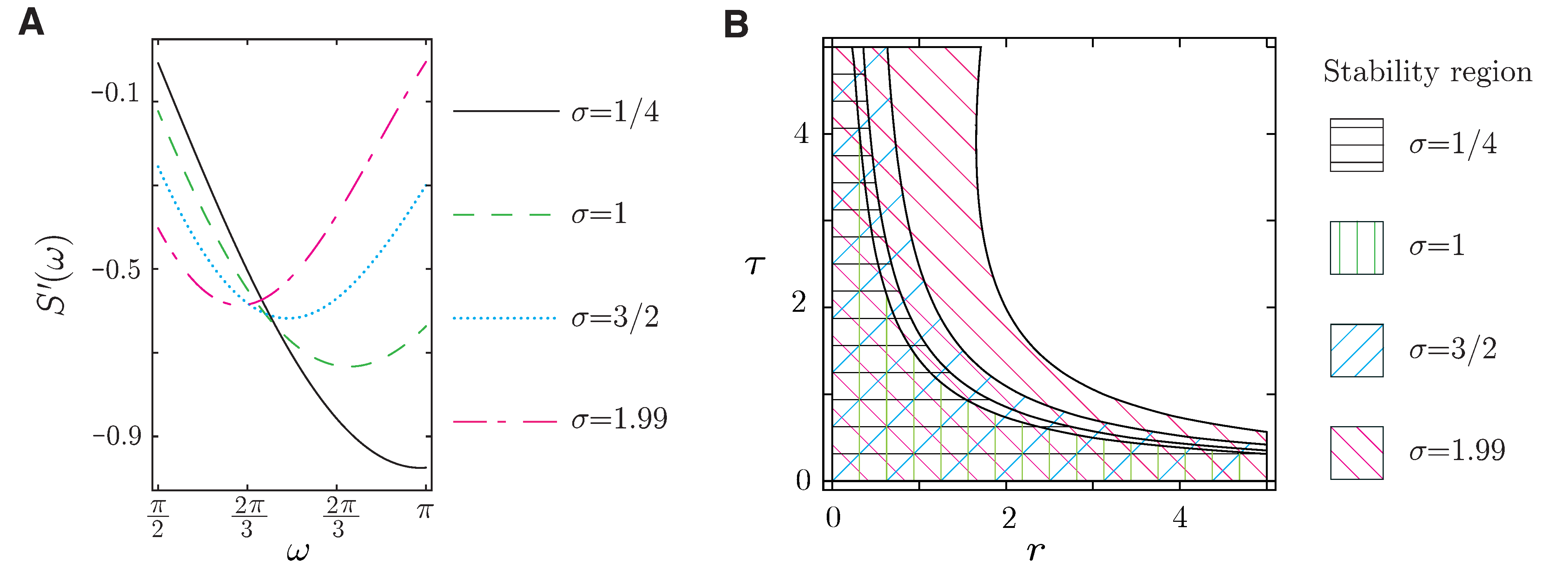}
          \caption{\textbf{Model \eqref{model} with uniform distribution kernel.}  
        \textbf{(A)}  The plot of $S'(\omega)$ given in \eqref{fun_s_prime}
        with different values of $\sigma\in(0,2)$.
        \textbf{(B)} The stability region of $n^*$ with different values of $\sigma\in(0,2)$ and fixing  the value of
other parameters $K=5$ and $D=3$.}
  \label{Fig:kernel_uniform_all_curves}
\end{figure*}

For $\omega\in(\pi/2,\pi)$
\[\ddc{\tau}{\omega} =-\cfrac{K}{r(n^*-K)} \left( \cfrac{ -\omega }{\sin^2(\omega)} + \cfrac{\cos(\omega)}{\sin(\omega)}   \right)>0\]
due to ${\cos(\omega)}/{\sin(\omega)}<0$. 
Moreover, we have
\begin{equation}\label{fun_s_prime}
    S'(\omega)=\cfrac{1}{\sigma  \omega ^2}
\Big(  {2 \sin \left(\sigma  \omega/2\right) (\omega  \cos (\omega )-\sin (\omega ))+\sigma  \omega  \sin (\omega ) \cos \left(\sigma  \omega/2\right)}\Big)<0
\end{equation}
when  $\omega\in(\pi/2,\pi)$ and $\sigma\in(0,2)$, see Fig. \ref{Fig:kernel_uniform_all_curves}A.

By Fixing  $K$ and  $D$, we plot the Hopf bifurcation curve in the $(r,\tau)$-plane using \eqref{condition_d_ge_1} with different values of $\sigma$ in Fig. \ref{Fig:kernel_uniform_all_curves}B.
We can see that  the  equilibrium $n^*$   is locally asymptotically stable  below the Hopf bifurcation curve and as $\tau$ increases when $r$ is fixed, $n^*$ becomes unstable above the Hopf bifurcation curve.
Furthermore, The figure shows that as $\sigma$ increases, the stability region (below  the Hopf bifurcation curve) increases. 

\begin{figure*}[hbt!] 
     \centering
         \includegraphics[width=1\textwidth]{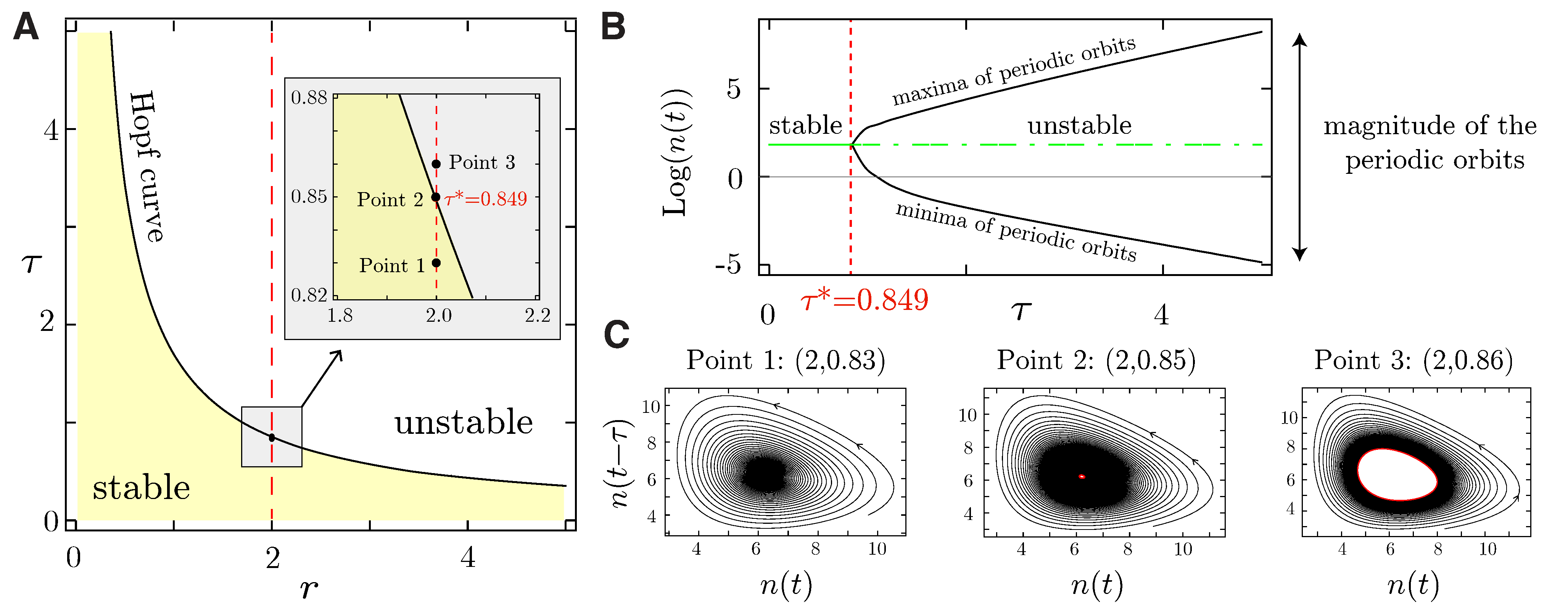}
          \caption{\textbf{Dynamics of model \eqref{model} with uniform distribution kernel.}  
        \textbf{(A)}  Stability region in $(r,\tau)$-plane with $\sigma=1$.
        \textbf{(B)} One-parameter bifurcation diagram when $r=2$ in (A). The positive equilibrium $n^*$ loses its stability at  $\tau^*=0.849$.
        \textbf{(C)} Phase portrait
 when $\tau=0.83,\,0.85$, and $0.86$  in (B). The value of
other parameters is $K=5$ and $D=3$.}
  \label{Fig:kernel_uniform}
\end{figure*}

For further discussion, we consider the case of $\sigma=1$ in Fig. \ref{Fig:kernel_uniform_all_curves} and study the dynamics of the model \eqref{model_uni} in Fig. \ref{Fig:kernel_uniform}. We can see that when fixing $r$ and increasing $\tau$, a limit cycle appears when crossing the Hopf bifurcation curve. Moreover, the magnitude of the limit cycle increases as $\tau$ increases.

 \subsection{Application 2: Dirac-Delta kernel}	
\label{sec_Dirac}

The Dirac-Delta kernel (Fig. \ref{Fig:kernel_forms}B) takes the form: 
    \begin{equation}\label{kernel_delta}
       g(s)=\delta_{\tau}(s)
       = \left\{\begin{array}{lll}0& \text{if}& s\ne \tau,\\ 
       \infty & \text{if} & s=\tau. \end{array}\right.
    \end{equation}
with mean time delay $\gavg=\tau$. 
When  $\tau=0$, model \eqref{model} reduces to an ordinary differential equation (ODE):
\begin{equation} \label{model_ode}
    \ddc{n(t)}{t}=r\, n(t)\, \left[   
1-\cfrac{n(t)}{K}
    \right]+D.
\end{equation}
While when $\tau>0$, model \eqref{model} reduces to a delay differential equation (DDE) with discrete time delay:
\begin{equation} \label{model_dde}
    \ddc{n(t)}{t}=r\, n(t)\, \left[   
1-\cfrac{n(t-\tau)}{K}
    \right]+D.
\end{equation}
Biologically, the distribution $\delta_{\tau}(s)$ means that the maximum influence on the  population density at present $t$ comes from a specific population density at last time $t-\tau$. 
Model \eqref{model_dde} with $D=0$ is studied in  \cite{beretta1987global,ruan2006delay}.

To study the stability with Dirc Delta kernel $\delta_{\tau}(s)$ defined  in \eqref{kernel_delta}, take $\sigma\to 0$, and hence, $C(\omega)=\cos(\omega)$ and $S(\omega)=\sin(\omega)$. Hence, 
when $D=0$, the equilibrium $n^*=K$  
is locally asymptotically stable if $\tau<\pi/(2r)$  and unstable when $\tau>\pi/(2r)$. Moreover, a Hopf bifurcation occurs at $\tau=\pi/(2r)$.
The result is consistent with \cite[Theorem 1]{ruan2006delay}.

 On the other hand, when $D>0$, there exists $\omega_0\in(\pi/2,\pi) $ such that $\omega_0={\rm arccos}(K/n^*-1)$, and  
the positive equilibrium $n^*$  
is locally asymptotically stable below the Hopf bifurcation curve defined by  $\tau=\omega_0 K/(\tau\, r\, n^*\,\sin(\omega_0))$. Moreover, $n^*$  is unstable above the Hopf bifurcation curve.

 \subsection{Application 3: Gamma distribution kernel}	
\label{sec_gamma}

The gamma distribution kernel (Fig. \ref{Fig:kernel_forms}C) can be written as:
\begin{equation}\label{kernel_gamma}
   ~\qquad   g(s)=\gammag(s):=\cfrac{\gamma^p s^{p-1} e^{-\gamma s}}{(p-1)!},\qquad  \gamma \ge0 \text{  and   }  p\in \mathbb{N}.
    \end{equation} 
The parameter  $p$ is the order of the delay kernel, and ${1}/{\gamma}$ is    the scale parameter. 
The mean time delay in this case is  $\gavg={p}/{\gamma}$.
When $g(s)=\gammag(s)$, model \eqref{model} reduces to an integro-differential equation (IDE) of the form
\begin{equation} \label{model_gamma}
    \ddc{n(t)}{t}=r\, n(t)\, \left[   
1-\cfrac{\gamma^p }{K\, (p-1)!}\int_{0}^{\infty} n(t-s)\,s^{p-1} e^{-\gamma s}\,\ds
    \right]+D.
\end{equation}
Using the linear chain trick \cite{macdonalds1978time} model \eqref{model_gamma} can be transformed to ODEs system of dimension  $p+1$ of the from 
\begin{align*}
  \ddc{n(t)}{t}&=r\, n(t)\, \left[   
1-\cfrac{x_p(t)}{K} 
    \right]+D.  \\
 \ddc{x_i(t)}{t}&= \gamma\, (x_{i-1}(t)-x_i(t)),\qquad i=1,2,\dots,p   
\end{align*}
where $x_0(t)=n(t)$ and
\[x_i(t)=\cfrac{\gamma^i }{(i-1)!}\int_{0}^{\infty} n(t-s)\,s^{i-1} e^{-\gamma s}\,\ds.\]
When $p=1$, the kernel $ g(s)$ is called \textit{exponential distribution} or \textit{weak delay kernel}. 
From a biological perspective, it shows that the maximum weighted response of population density  comes from the present population density. 
While $ g(s)$ is called \textit{strong delay kernel} when $p=2$. 
Biologically, it means that the maximum influence on the population density at any time $t$ is determined by the density of the population at the preceding time $t- 1/\gamma$. See Fig. \ref{Fig:kernel_forms}C.

The normalized gamma distribution has the form \[\widehat{g}(s)=\cfrac{p}{\gamma}\,\gammag(p s/\gamma)=\cfrac{p^p s^{p-1}\e^{ps}}{(p-1)!}.\]
Consequently, the characteristic equation is
 \begin{equation}\label{characteristic_eq_gamma}
    \Delta(\lambda)=\lambda+\cfrac{p\, r}{\gamma}  \left( \cfrac{n^*-K}{K}\right)+ 
    \cfrac{p\, r\, n^*}{\gamma\, K} \left( \cfrac{p}{\lambda+p}  \right)^p=0.
\end{equation}
Following \cite{campbell2009approximating} we have
\begin{neweq_non}
C(\omega)&= \re \left[ \frac{p^p}{(p-1)!} \int_0^{\infty}  s^{p-1}\e^{-(p+i\omega) s} \ds  \right]=
\left( 1+\cfrac{\omega^2}{p^2}  \right)^{-p} \re \left( 1-\cfrac{i\omega}{p}     \right)^p\\ 
&= \left( 1+\cfrac{\omega^2}{p^2}  \right)^{-p} \sum _{j=0}^{\lfloor \frac{p}{2}\rfloor } \binom{p}{2j} \frac{(-1)^{j}{\omega }^{2j}}{{p}^{2j}}
\end{neweq_non}
and
\begin{neweq_non}
S(\omega)&=
\left( 1+\cfrac{\omega^2}{p^2}  \right)^{-p} \im \left( 1-\cfrac{i\omega}{p}     \right)^p\\
&=\left( 1+\cfrac{\omega^2}{p^2}  \right)^{-p} \sum _{j=0}^{\lfloor \frac{p-1}{2}\rfloor } \binom{p}{2j+1} \frac{(-1)^{j}{\omega }^{2j+1}}{{p}^{2j+1}}.
\end{neweq_non}
For $p=1$, we obtain
\begin{neweq_non}
C(\omega)=\cfrac{1}{1+\omega^2}  \qquad\text{and}\qquad S(\omega)=\cfrac{\omega}{1+\omega^2}.
\end{neweq_non}
When $p=2$, we have
\begin{neweq_non}
C(\omega)=\cfrac{1-\omega^2/4}{(1+\omega^2/4)^2} \qquad\text{and}\qquad S(\omega)=\cfrac{\omega}{(1+\omega^2/4)^2},
\end{neweq_non}
and $p=3$ gives
\begin{neweq_non}
C(\omega)=\cfrac{1-\omega^2/3}{(1+\omega^2/9)^3} 
\qquad\text{and}\qquad 
S(\omega)=\cfrac{\omega(1-\omega^2/27)}{(1+\omega^2/9)^3}.
\end{neweq_non}

\subsubsection[The case of D=0]{The case of $D=0$.}

The case of $p=1$ is straightforward because  $C(\omega)=0$  has no real roots.  Thus, the positive equilibrium $n^*=K$ is always locally asymptotically stable.

Recall that the mean time delay is  $\gavg={p}/{\gamma}$.  
When $p=2$, $C(\omega)=0$ at $\omega_0=2$. Thus, $n^*=K$ is locally asymptotically stable if 
\[\cfrac{2r}{\gamma}<\cfrac{2}{S(2)}=4
\qquad\Rightarrow\qquad
r<2\, \gamma.\]
Similarly, when $p=3$, we obtain that $\omega_0=\sqrt{3}$ and $n^*=K$ is locally asymptotically stable when \[r<\cfrac{8}{9}\, \gamma.\] 
It is easy to check that $C'(\omega_0)<0$ and a Hopf bifurcation occurs when crossing the imaginary axis.
These results are consistent with the results in \cite{sawada2022stability} where the stability condition is 
      \[r< \cfrac{\tan\left({\pi}/{(2p)} \right)}{\cos^p\left({\pi}/{(2p)} \right)}\,\,\gamma.\].

\begin{figure*}[hbt!]
     \centering
         \includegraphics[width=1\textwidth]{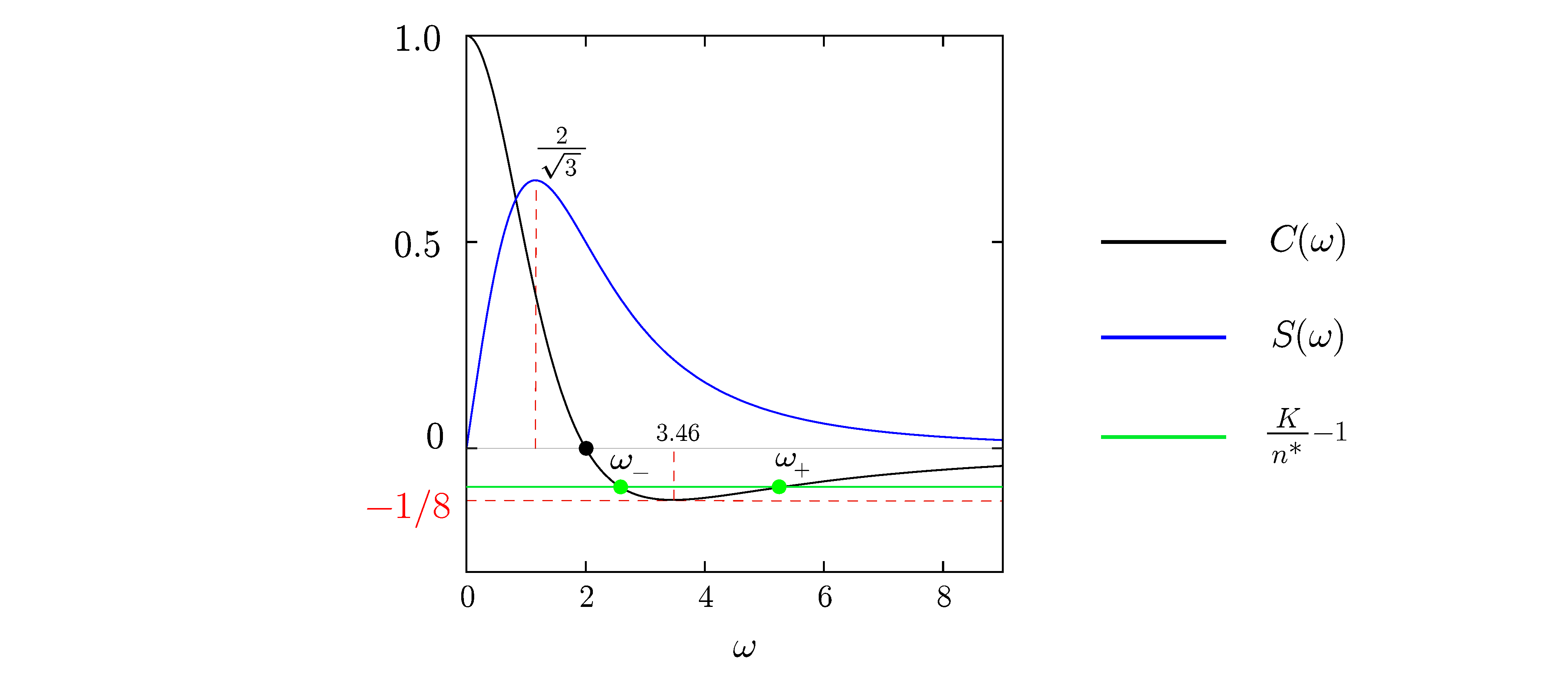}
          \caption{\textbf{Existence of pure imaginary roots $\omega_{\pm}$.}  
           The intersection of the curve $C(\omega)$ and horizontal line $K/n^*-1$ when $r>r^*$ defined in \eqref{condition_r}.}
  \label{Fig5_gamma_CS}
\end{figure*}


\subsubsection[The case of D>0]{The case of $D>0$.}

The case of $p=1$ is straightforward because  $C(\omega)=K/n^*-1<0$ has no real  roots.  Thus, the positive equilibrium $n^*$ is always locally asymptotically stable. 

When $p=2$, $-1/8=C(3.46)\le C(\omega)\le 1$, and hence, the equation $C(\omega)=K/n^*-1$ has a solution if
\begin{equation}\label{condition_r}
    r>\cfrac{49\, D}{8\, K}:=r^* 
\end{equation}
due to $K/n^*-1<0$.

When the condition \eqref{condition_r} holds, the curve of $C(\omega)$ intersects the horizontal line $K/n^*-1$ at two points
\[\omega_{\pm}=\sqrt{\frac{2 n^*}{n^*-K}\pm2 \sqrt{\frac{n^* (8 K-7 n^*)}{(n^*-K)^2}}-4},\]
where $2<\omega_-<3.46<\omega_+$ because $C(\omega)<0$ for $\omega>2$,  see Fig. \ref{Fig5_gamma_CS}. It is easy to check that 
\[S'(w)=\frac{16(4-3 w^2)}{\left(w^2+4\right)^3}<0\qquad\text{for}\qquad w>\cfrac{2}{\sqrt{3}}.\]
Thus, $S'(\omega_{\pm})<0$. Moreover, we have 
\begin{neweq_non}
\ddc{\gavg}{\omega}=\cfrac{16\,n^*\, \omega}{\gamma\,(n^*-K)\left(w^2+4\right)^2}>0\quad\text{for}\quad w>0.
\end{neweq_non}
Consequently, if \eqref{condition_d_ge_1} holds, then the equations
\[\cfrac{1-\omega^2/4}{(1+\omega^2/4)^2}=\cfrac{K}{n^*}-1
\qquad\qquad \text{and}\qquad\qquad
\cfrac{1}{(1+\omega^2/4)^2}=\cfrac{\gamma\, K}{2 r\, n^*}
\]
define parametric equations for Hopf bifurcation curve  with $\omega>0$.

\begin{figure*}[hbt!]
     \centering
         \includegraphics[width=1\textwidth]{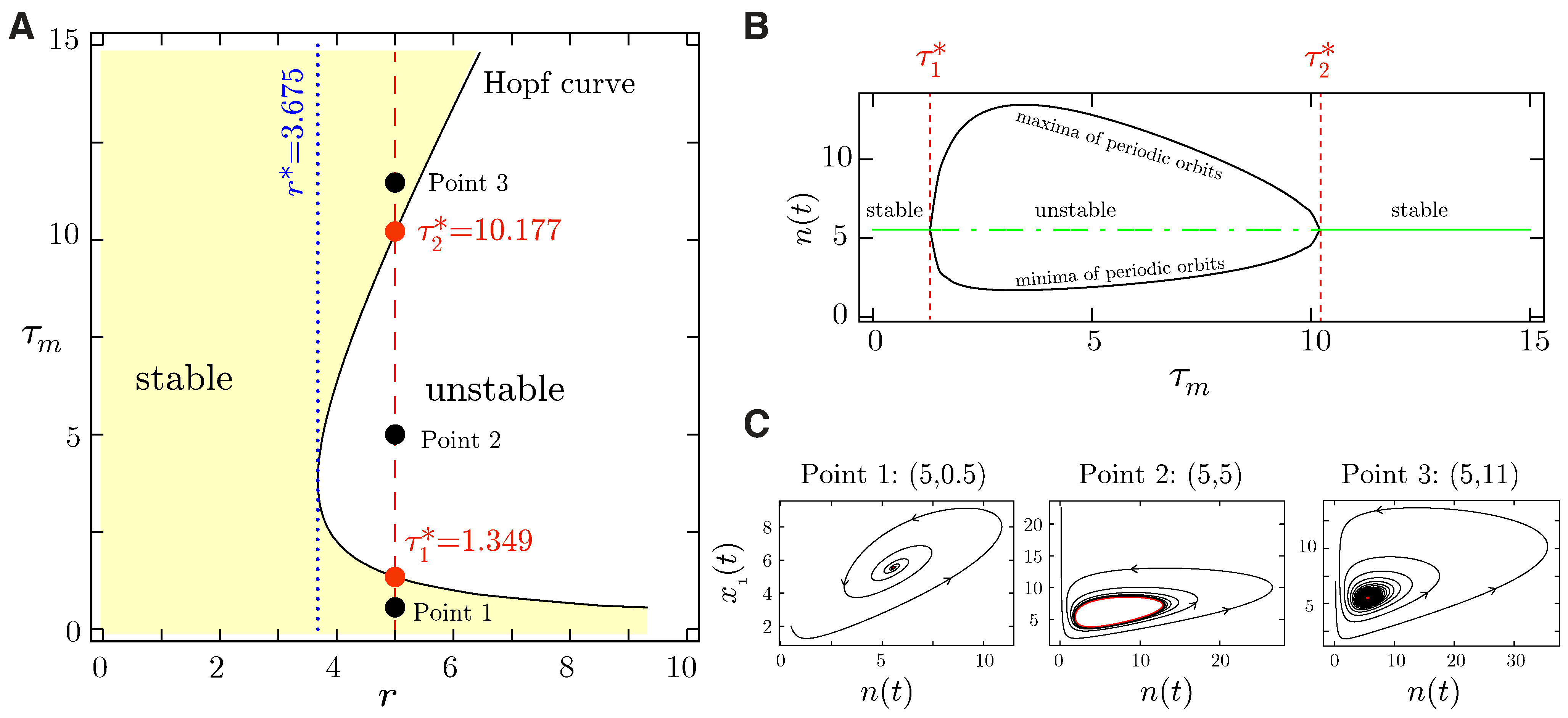}
          \caption{\textbf{Dynamics of model \eqref{model} with gamma distribution kernel and $p=2$.}  
        \textbf{(A)}  Stability region in $(r,\gavg)$-plane.      
        \textbf{(B)} One-parameter bifurcation diagram when $r=5$ in (A). The positive equilibrium $n^*$ loses its stability at  $\tau_1^*=1.349$ and gains it again at $\tau_2^*=10.177$. 
          \textbf{(C)} Phase portrait
 when $\tau=0.5,\,5$, and $11$  in (B). 
 The value of
other parameters is $K=5$ and $D=3$.}
  \label{Fig:kernel_gamma}
\end{figure*}

 In Fig. \ref{Fig:kernel_gamma}, we plot the Hopf  bifurcation curve in  $(r,\gavg)$-plane by fixing $K$ and $D$. For $r<r^*$,  the  equilibrium $n^*$ is stable for all $\gavg$. On the other hand, 
    when  $r>r^*$ and there exist  $0<\tau_1^*<\tau_2^*$  such that 
    as $\gavg$ increases (or $\gamma$ decreases) the  equilibrium $n^*$ is stable for $\gavg<\tau_1^*$,  loses the
stability when $\gavg\in(\tau_1^*,\tau_2^*)$, and back to stable for $\gavg>\tau_2^*$. 

When the order of the delay kernel is $p=3$, we have  $-1/4=C(3)\le C(\omega)\le 1$. Recall that $K/n^*-1<0$. Thus, the equation $C(\omega)=K/n^*-1$ has a solution if
\begin{equation}\label{condition_r_p_3_L}
    r>\cfrac{9\, D}{4\, K}:=\underline{r}.
\end{equation}
Consequently, when condition \eqref{condition_r_p_3_L} holds, the solution of the equation $C(\omega)=K/n^*-1$ has two solutions $\omega_{\pm}$ such that $\sqrt{3}<\omega_-<3<\omega_+$.

It is easy to check that
\[S'(\omega)=\cfrac{81 \left(w^4-54 w^2+81\right)}{\left(w^2+9\right)^4}\]
and
\begin{neweq_non}
\ddc{\gavg}{\omega}=-\cfrac{2^4\,3^7\,n^*\, \omega}{\gamma\,(n^*-K)\left(w^2-27\right)\left(w^2+9\right)^2}.
\end{neweq_non}

When $\omega\to 3\sqrt{3}^-$, we have $\tau\to\infty$ and 
\begin{equation}\label{condition_r_p_3_U}
r\to\cfrac{49\, D}{8\, K}:=\overline{r}.
\end{equation}
Thus, when fixing $K$ and $D$,  the  Hopf  bifurcation curve in  $(r,\gavg)$-plane has a vertical asymptote at $r=\overline{r}$.

\begin{figure*}[hbt!]
     \centering
         \includegraphics[width=1\textwidth]{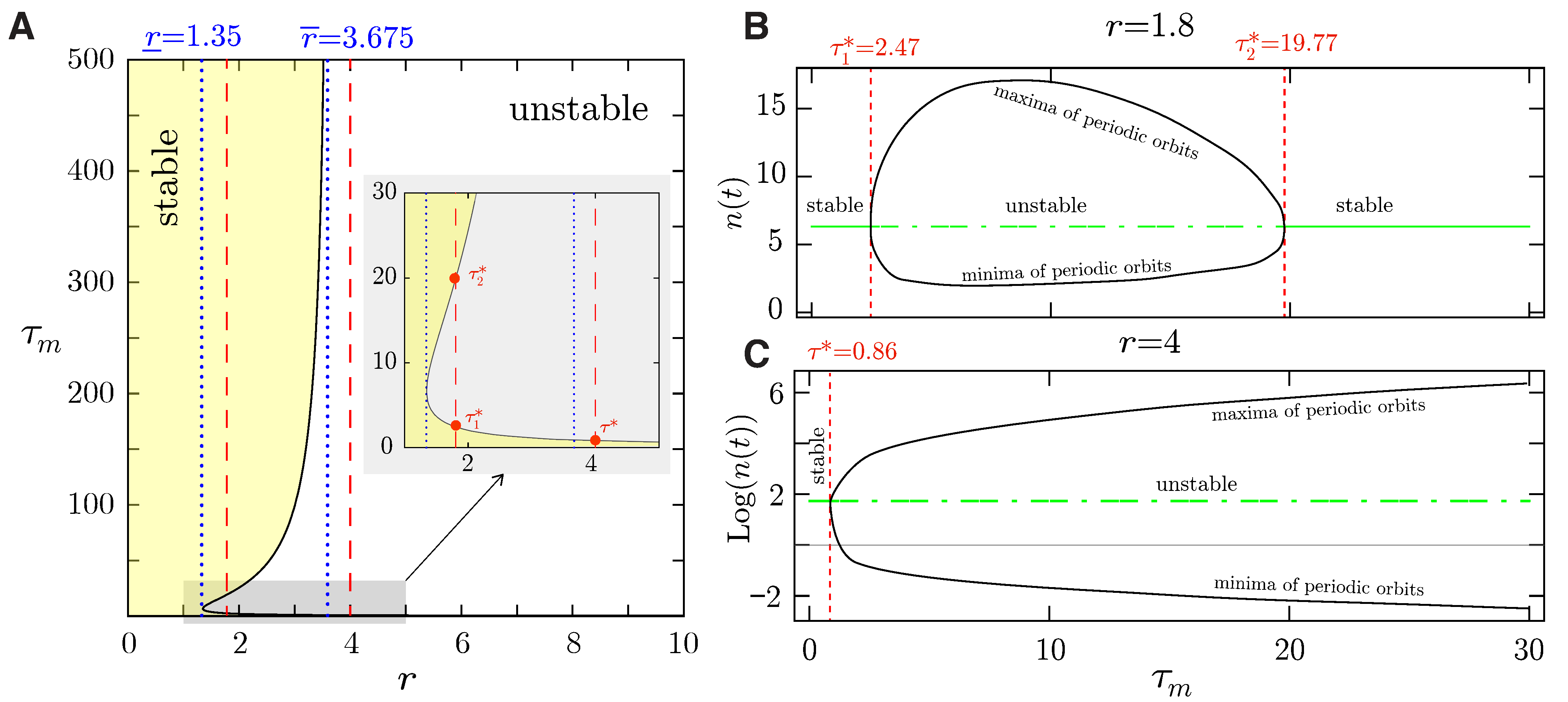}
          \caption{\textbf{Dynamics of model \eqref{model} with gamma distribution kernel and $p=3$.}  
      \textbf{(A)}  Stability region in $(r,\gavg)$-plane.      
        \textbf{(B)} One-parameter bifurcation diagram when $r=1.8$ in (A). The positive equilibrium $n^*$ loses its stability at  $\tau_1^*=2.46$ and gains it again at $\tau_2^*=19.77$. 
         \textbf{(C)} One-parameter bifurcation diagram when $r=4$ in (A). The positive equilibrium $n^*$ loses its stability at  $\tau^*=0.86$. 
 The value of
other parameters is $K=5$ and $D=3$.}
  \label{Fig:kernel_gamma_p_3}
\end{figure*}

In Fig. \ref{Fig:kernel_gamma_p_3}, we plot the Hopf  bifurcation curve in  $(r,\gavg)$-plane by fixing $K$ and $D$. For $r<\underline{r}$,  the  equilibrium $n^*$ is stable for all $\gavg$. On the other hand, 
    when  $\underline{r}<r<\overline{r}$ and there exist  $0<\tau_1^*<\tau_2^*$  such that 
    as $\gavg$ increases (or $\gamma$ decreases) the  equilibrium $n^*$ is stable for $\gavg<\tau_1^*$,  loses the
stability when $\gavg\in(\tau_1^*,\tau_2^*)$, and back to stable for $\gavg>\tau_2^*$. Moreover, we see that when $r>\overline{r}$, there exists $\tau^*$ such that the  equilibrium $n^*$   is locally asymptotically stable for $\gavg<\tau^*$ and unstable for  $\gavg>\tau^*$.

\section{Conclusions}
\label{sec_conclusions}

In this paper, we have studied the stability of a time-delayed single-species logistic model with a general distribution delay kernel and an inflow of nutritional resources at a constant rate $D\ge 0$.
By normalizing the time delay, we have studied the stability of the positive equilibrium point and provided precise conditions for the linear stability and the occurrence of Hopf bifurcation.
To study the influence of population density at any previous time on its density at present, we have applied the general results to three delay distribution kernels: Uniform, Dirac-delta, and gamma, where each distribution has its distinctive biological interpretation.

In the case of zero delay, we have proved that the positive equilibrium is always stable with $D\ge 0$.
For Uniform and Dirac-delta distributions, we have found the dynamics are similar when $D=0$ or $D>0$, where  the positive equilibrium is stable for a relatively small delay and then loses its stability through the Hopf bifurcation when the mean delay increases.
On the other hand, we have shown that the dynamics are different for the gamma distribution when $D>0$, which is affected by the delay order. 
We have proved that the positive equilibrium is always stable when the delay order is one (weak gamma distribution or exponential distribution).
When the order of delay is two or three, we have shown that there is stability switching of the positive equilibrium resulting from the increase of the value of $\gavg$, in the sense that the positive equilibrium is stable for a relatively short period. Then, it loses stability via Hopf bifurcation as $\gavg$ increases. After then, it stabilizes again with an increase in $\gavg$. 
The main difference  is that for relatively large $\gavg$ and intrinsic growth rate, the positive equilibrium can be stable when delay order is two, but it will be unstable when delay order is three.
We have found that the results are consistent with the parts in the literature \cite{sawada2022stability,ruan2006delay}.

For future work, we will consider a multi-patch time delayed logistic equation with migration \cite{al2016prey} where each patch follows a different distributed delay kernel. Then, study the influence of delay in the density of population in each patch in the presence of migration.

\bibliographystyle{ieeetr}   
\bibliography{References.bib}

\end{document}